%% file: Seifert_accepted.tex
\documentclass{amsart}
\input{command}

\begin{document}

\title{The Guts of nearly fibered knots}

%    Remove any unused author tags.

%    author one information
\author{Zhenkun Li}
\address{Department of Mathematics, University of South Florida}
\curraddr{}
\email{zhenkun@usf.edu}
\thanks{}

%author two information
\author{Fan Ye}
\address{Department of Mathematics, Harvard University}
\curraddr{}
\email{fanye@math.harvard.edu}
\thanks{}

\keywords{}
\date{}
\dedicatory{}
\begin{abstract}
The guts of a knot is an invariant defined for the knot complement by Agol-Zhang. Nearly fibered knots, which are defined as knots whose Floer homology has dimension two in the top Alexander grading, were introduced by Baldwin-Sivek. In this note, we provide three models for the guts of nearly fibered knots in the $3$-sphere. As a corollary, the nearly fibered condition can be purely topologically characterized and is independent of the specific version of Floer theory.
\end{abstract}

\maketitle

%\tableofcontents%table of contents
% \newpage

%————Start from here————
\section{Introduction}

By work of Ghiggini \cite{ghiggini2008knot}, Ni \cite{ni2007knot} and Kronheimer-Mrowka \cite{kronheimer2010knots}, a knot $K\subset S^3$ is fibered if and only if its knot homology in any branch of Floer theory is $1$-dimensional in the top Alexander grading. Hence it is natural to ask what happens if the top grading summand of the knot homology is $2$-dimensional. Recently, Baldwin-Sivek in \cite{baldwin2022nearly} introduced the following definition.

\begin{defn}
	A knot $K\subset S^3$ is said to be {\bf nearly fibered} (in the Heegaard Floer sense) if 
	$$\widehat{HFK}(S^3,K,g(K);\mathbb{Q})\cong\mathbb{Q}^2.$$
\end{defn}

Their definition is stated with Heegaard Floer theory, but we can also define nearly fibered knots in the instanton sense by requiring
$$KHI(S^3,K,g(K))\cong \mathbb{C}^2,$$
where $KHI$ denotes the instanton knot homology \cite{kronheimer2010knots} of $K\subset S^3$.

In this note, we show that the nearly fibered condition has a purely topological characterization, and is independent of the branches of Floer theory. To better describe this criterion, we use the notion of guts of knots recently introduced by Agol-Zhang \cite{Agol2022guts}. 

Given a knot $K\subset S^3$, we can view its complement $S^3\backslash N(K)$ as a sutured manifold with its whole boundary being the suture. We can pick a maximal collection of pair-wise disjoint and pair-wise non-parallel minimal-genus Seifert surfaces $S$ of $K$, and perform a sutured manifold decomposition
\begin{equation}\label{eq: decompose along Seifert surface}
    S^3\backslash N(K)\stackrel{S}{\leadsto}(M^\p,\ga^\p).
\end{equation}
We can then pick a maximal collection of pair-wise disjoint and pair-wise non-parallel non-trivial product annuli $A$ inside $(M^\p,\ga^\p)$ and perform a second sutured manifold decomposition\footnote{Here, note that in Agol and Zhang's paper \cite{Agol2022guts}, they also require to decompose along non-trivial product disks to obtain the guts. In this paper we drop the step of decomposing along possible product disks because of \cite[Lemma 2.13]{juhasz2010polytope}: the only two taut balanced sutured manifolds that admit no non-trivial product annulus but admit non-trivial product disks are both product sutured manifolds, and hence are actually the components to be dropped when obtaining guts.}
\begin{equation}\label{eq: decompose along product annuli}
    (M^\p,\ga^\p)\stackrel{A}{\leadsto}(M,\ga)\sqcup (M_1,\ga_1).
\end{equation}
Here $(M_1,\ga_1)$ is a product sutured manifold, and no components of $(M,\ga)$ is a product.
\begin{defn}[\cite{Agol2022guts}]
    The guts of a knot $K\subset S^3$ is defined to be the sutured manifold $(M,\ga)$.
\end{defn}

\bthm[{\cite[Theorem 1.1]{Agol2022guts}}]\label{thm: well-defined of guts}
The guts of a knot $K\subset S^3$ is well-defined, {\it i.e.}, independent of the choices of maximal collections of Seifert surfaces and product annuli in the construction.
\ethm

% Roughly speaking, the \textbf{guts} of a knot $K\subset S^3$ is obtained from the knot complement $S^3\backslash N(K)$ by a two-step sutured manifold decompositions: we first decompose $S^3\backslash N(K)$ along a maximal collection of pair-wise disjoint minimal-genus Seifert surfaces and then decompose the resulting sutured manifold along a maximal collection of pair-wise disjoint product annuli. Then the non-product pieces are defined to be the guts. For the definition of sutured manifold decompositions, readers are referred to Gabai \cite{gabai1983foliations}. 

In this note, we prove that a knot $K\subset S^3$ is nearly fibered if and only if its guts falls into one of the three basic models described below. Note that in our paper we only state and prove the theorem in instanton theory, but a similar argument applies to Heegaard Floer theory as well.

\bthm\label{thm: rank-2 top grading}
Suppose $K\subset S^3$ is a knot of genus $g$. Let $(M,\ga)$ be its guts. Then we have
$$KHI(S^3,K,g)\cong\mathbb{C}^2,$$
if and only if its guts $(M,\ga)$ falls into one of the following three models up to orientation reversal of the ambient $3$-manifold:%\footnote{One might worry that the orientation of the suture also matters. Though in all three models, $\partial M$ is a torus, so $\gamma$ and $-\gamma$ are isotopic on $\partial M$.}

\begin{enumerate}
	\item [(M1)] $M$ is a solid torus and $\ga$ consists of four longitudes.
		\item [(M2)] $M$ is a solid torus and $\ga$ consists of two curves of slope $2$.
		\item [(M3)] $M$ is the complement of the right handed trefoil and $\ga$ consists of two curves of slope $2$.
\end{enumerate}
\ethm

We have the following corollary.

\bcor
A knot is nearly fibered in the instanton sense if and only if it is nearly fibered in the Heegaard Floer sense.
\ecor

{\bf Acknowledgement} The authors would like to thank John A. Baldwin and Steven Sivek for helpful conversations and for telling us about their work on classifying the genus-one nearly fibered knots.

\section{Proofs and comments}\label{sec: rank-2 top grading}
% In this section, we first define the guts of a knot and then prove the main theorem. Notations for instanton Floer homology we use in this paper can be found in any of our previous papers, for example, \cite{LY2020}.

\bpf[Proof of Theorem \ref{thm: rank-2 top grading}]
We first prove the necessary condition. Suppose $K\subset S^3$ is a genus-$g$ nearly fibered knot. Let $(M,\ga)$ be its guts. We first study the sutured manifold decomposition in (\ref{eq: decompose along Seifert surface}).

\begin{claim}
	Any maximal collection of pair-wise disjoint and pair-wise non-parallel minimal-genus Seifert surfaces in fact contains only one Seifert surface.
\end{claim}

\begin{proof}[Proof of Claim 1]Suppose $S$ is a minimal-genus Seifert surface of $K$. We can perform a sutured manifold decomposition of $S^3\backslash N(K)$ along $S$:
$$S^3\backslash N(K)\stackrel{S}{\leadsto}(S^3\backslash [-1,1]\times S,\{0\}\times\partial S).$$
By the proof of \cite[Proposition 7.16]{kronheimer2010knots}, we know that there is an isomorphism
\begin{equation}\label{eq: top grading and Seifert surface complement}
    SHI(S^3\backslash [-1,1]\times S,\{0\}\times\partial S)\cong KHI(S^3,K,g)\cong\mathbb{C}^2.
\end{equation}
If there is another minimal-genus Seifert surface $S^\p$ that is disjoint from $S$ and is not parallel to $S$, then $S^\p$ also induces a non-boundary parallel surface in $(S^3\backslash [-1,1]\times S,\{0\}\times\partial S)$, which implies the sutured manifold is not horizontally prime. From the instanton version of \cite[Proposition 6.5 and Proposition 6.6]{kronheimer2010knots}, we know that one of the two pieces obtained from $(S^3\backslash [-1,1]\times S,\{0\}\times\partial S)$ by cutting along $S^\p$ must have $1$-dimensional sutured homology, because $2$ is a prime number. From \cite[Theorem 7.18]{kronheimer2010knots}, that piece is a product sutured manifold, which contradicts the assumption that $S^\p$ is not parallel to $S$.
\end{proof}

Now Claim 1 above and Theorem \ref{thm: well-defined of guts} imply that the sutured manifold $(M^\p,\ga^\p)$ in (\ref{eq: decompose along Seifert surface}) can be taken to be simply the complement of $S$:
$$(M^\p,\ga^\p)=(S^3\backslash [-1,1]\times S,\{0\}\times\partial S).$$

Next, we study the sutured manifold decomposition (\ref{eq: decompose along product annuli}). Note that by construction $(M_1,\ga_1)$ is a product sutured manifold, so from \cite[Theorem 7.18]{kronheimer2010knots} and the instanton version of \cite[Proposition 6.5 and Proposition 6.7]{kronheimer2010knots}, we know that 
$$SHI(M,\ga)\cong SHI(S^3\backslash [-1,1]\times S,\{0\}\times\partial S)\cong \mathbb{C}^2.$$
Also, the same argument as above shows that $(M,\ga)$ is horizontally prime. Thus we conclude that $(M,\ga)$ is reduced in the sense of \cite[Definition 2.12]{juhasz2010polytope}. Then \cite[Corollary 1.16]{li2019decomposition} applies and we conclude that
$$b_1(M)=b^1(M)\leq 2-1=1.$$
We claim that $g(\partial M) = b_1(M)$. Indeed, we know that $(M,\gamma)$ is obtained from the knot complement by decomposition. So \cite[Lemma 5.1]{juhasz2010polytope} implies that $H_2(M) = 0$. As a result, by the universal coefficient theorem and the Poincar\'{e} duality, we have
\[
H_1(M,\partial M;\mathbb{Q})\cong H^2(M;\mathbb{Q}) \cong H_2(M;\mathbb{Q}) = 0.
\]
Hence the long exact sequence of the pair $(M,\partial M)$ implies that the map
\[
i_*: H_1(\partial M;\mathbb{Q}) \to H_1(M;\mathbb{Q})
\]
is surjective. Hence the `half lives and half dies' theorem in $3$-dimensional topology implies that
\[
g(\partial M) = \frac{1}{2}b_1(\partial M) = b_1(M).
\]

Now $g(\partial M) = b_1(M) \leq 1$. If $g(\partial M)=0$, since $M$ is irreducible, we know $M=B^3$. Then for any possible $\ga$ on $\partial M$, we cannot have $SHI(M,\ga)\cong\mathbb{C}^2$.
Hence we must have $\partial M\cong T^2$. It is well-known that any (smooth) torus in $S^3$ bounds a solid torus. Hence we have two cases.

\vspace{1em}\noindent%
{\bf Case 1.} The manifold $M$ is a solid torus. 

The instanton Floer homology of any sutured solid torus can be found in \cite[Section 4.3]{li2019direct}. So the only two models are the ones as stated in (M1) and (M2).

\vspace{1em}\noindent%
{\bf Case 2.} The manifold $S^3\backslash M$ is a solid torus, {\it i.e.}, there is a knot $J\subset S^3$ so that $M\cong S^3\backslash N(J)$. 

Suppose $\ga$ has $2n$ components. Let $\ga_2$ be the union of two adjacent components of $\ga$, which are necessarily oppositely oriented. Next, we make the following claim.

\begin{claim}
	Suppose $(M,\gamma)$ is a balanced sutured manifold and assume that three components of $\gamma$ are parallel disregarding the orientation. Write $\gamma_3$ to be the disjoint union of these three copies. Note that two components of $\gamma_3$ are coherently oriented and are opposite to the third. Let $\gamma_1$ be either of the two coherently oriented components and write $\gamma' = (\gamma \backslash \gamma_3) \cup \gamma_1$. Then we have
\[
SHI(M,\gamma) = SHI(M,\gamma')\otimes \mathbb{C}^2.
\]

\end{claim}
\bpf[Proof of Claim 2] The proof essentially follows from the proof of \cite[Theorem 3.1]{kronheimer2010instanton}. There exists an embedded annulus $A\subset \partial M$ such that $A$ contains $\gamma_3$ and each component of $\gamma_3$ is a core of $A$. Push the interior of $A$ into the interior of $M$ to produce a properly embedded annulus. Fix any orientation of $A$. Then there is a product annulus decomposition
\[
(M,\gamma)\stackrel{A}{\leadsto}(V,\gamma^4)\sqcup (M,\gamma'),
\]
where $V$ is a solid torus and $\gamma^4$ consists of four longitudes (there is a unique way, up to isotopy, to make $(V,\gamma^4)$ a balanced sutured manifold). Now an instanton version of \cite[Proposition 6.7]{kronheimer2010knots} implies that
\[
SHI(M,\gamma) = SHI(M,\gamma')\otimes SHI(V,\gamma^4).
\]
In the proof of \cite[Theorem 3.1]{kronheimer2010instanton}, Kronheimer and Mrowka already computed that
\[SHI(V,\gamma^4) \cong \mathbb{C}^2\]
and hence we are done.
\epf

Applying Claim 2 repetitively, we conclude that
$$SHI(M,\ga)\cong \mathbb{C}^{2^{n-1}}\otimes SHI(M,\ga_2).$$
Since
$$SHI(M,\ga)\cong \mathbb{C}^2,$$
either $n=2$ and $SHI(M,\ga_2)\cong \mathbb{C}$, or $n=1$. For the former case, from \cite[Theorem 7.18]{kronheimer2010knots} we know $M$ must also be a solid torus which reduces to Case 1. For the latter case, we further divide it into two sub-cases.

\vspace{1em}\noindent%
{\bf Case 2.1.} Each component of $\ga$ represents a generator of $\ker i_*\subset H_1(\partial M)\cong \intg^2$, where 
$$i_*: H_1(\partial M)\ra H_1(M)$$
is the map induced by the natural inclusion
$$i: \partial M\hookrightarrow M.$$

In this case, first recall that $M$ is a knot complement $S^3\backslash N(J)$ and hence $H_2(M,\partial M)$ is generated by a minimal-genus Seifert surface $T$ of the knot $J$. The assumption of Case 2.1 is equivalent to that $\gamma$ is parallel to $\partial T\subset \partial M$. We can assume that $\partial T\cap \gamma = \emptyset $. If $T$ is a disk, then $R(\ga)$ is compressible and $SHI(M,\ga)=0$ by the adjunction inequality ({\it cf.} \cite[Proposition 7.5]{kronheimer2010knots}). From now on we assume that $T$ has genus at least $1$. We know from \cite[Lemma 6.2]{li2019decomposition} that we have two taut decompositions
$$(M,\ga)\stackrel{\pm T}{\leadsto}(M_{\pm},\ga_{\pm}).$$

We make the following claim.

\begin{claim}
	We have an inclusion
\[
SHI(M_+,\ga_+)\oplus SHI(M_-,\gamma_-)\hookrightarrow SHI(M,\ga).
\]

\end{claim}
\bpf[Proof of Claim 3] We adopt the idea in \cite[Section 3]{li2019direct}. We isotope $T$ to $T^{\pm}$ such that the decomposition of $(M,\gamma)$ along $T^+$ is $(M_+,\ga_+)$ and the decomposition of $(M,\gamma)$ along $-T^-$ is $(M_-,\ga_-)$. $T^{\pm}$ are called positive and negative stabilizations of $T$ as in \cite[Definition 3.1]{li2019direct}, and we know that $-(T^-) = T^+$. By \cite[Theorem 3.4]{li2019direct}, each $T^{\pm}$ induces a $\intg$-grading on $SHI(M,\gamma)$. Then \cite[Lemma 4.2]{li2019direct} implies that
\[
SHI(M,\gamma, T^+, g(T)) \cong SHI(M_+,\gamma_+)
\]
and
\[
\begin{aligned}
	SHI(M,\gamma, T^-, -g(T)) &= SHI(M,\gamma, -T^-, g(T))\\
	&\cong SHI(M_-,\gamma_-).
\end{aligned}
\]
\cite[Proposition 4.1]{li2019direct} implies that\footnote{Note when reversing the orientation of the manifold and the suture, positive and negative stabilizations of $T$ are also switched.}
\[
SHI(M,\gamma, T^-, -g(T)) \subset SHI(M,\gamma, T^+, 1-g(T)).
\]
Hence we are done since $g(T)\neq 1-g(T)$.
\epf

Observe that both $\ga_{+}$ and $\ga_-$ contain at least three components that are parallel to each other. Let $\gamma_{\pm}'$ be the suture obtained from $\gamma_{\pm}$ by replacing three copies with one copy. Applying Claim 2, we know that
\[
SHI(M_{\pm},\gamma_{\pm})\cong SHI(M_{\pm},\gamma_{\pm}')\otimes \mathbb{C}^2.
\]
Tautness together with \cite[Theorem 7.12]{kronheimer2010knots} then implies
$$\dim SHI(M_{\pm},\ga_{\pm})\geq 2.$$
As a result, we have
$$\dim SHI(M,\ga)\geq 4,$$
which leads to a contradiction in this case.

\vspace{1em}\noindent%
{\bf Case 2.2.} Components of $\gamma$ do not represent generators of $\ker i_*$. 

Let $Y$ be the Dehn filling of $M$ along a component of $\ga$. We make the following claim.

\begin{claim}
	We have $\dim I^{\sharp}(Y) = 2$.
\end{claim}
\bpf[Proof of Claim 4] 

In order to prove Claim 4, we need the following three facts.

\begin{enumerate}
	\item [(1)] We have $\dim I^{\sharp}(Y)\neq 0$.
	\item [(2)] We have $\dim I^{\sharp}(Y) = \dim SHI(Y(1))\leq \dim SHI(M,\ga)=2$.
	\item [(3)] We have $\dim I^{\sharp}(Y)\equiv \dim SHI(M,\ga)~{\rm mod}~2$.
\end{enumerate}

To show (1), note that the fact $[\gamma]\notin \ker i_*$ implies that $Y$ is a rational homology sphere. Hence by \cite[Corollary 1.4]{scaduto2015instanton}, we know $\dim I^{\sharp}(Y)\neq 0$. 

To show (2), recall that $M= S^3\backslash N(J)$ is the knot complement and $\gamma$ has two components. Let $\gamma_0\subset \gamma$ be any component. We can attach a $3$-dimensional $2$-handle along $\gamma_0$. The resulting manifold is $Y\backslash B^3$. Hence we have a balanced sutured manifold
\[Y(1) = (Y\backslash B^3,\gamma \backslash \gamma_0).\]
Now let $T$ be the cocore arc of the $2$-handle. This arc $T$ is a vertical tangle inside $Y(1)$ as in \cite[Definition 1.1]{xie2019tangle}. Now observe that $(M,\gamma)$ can be obtained from $Y(1)$ by removing $T$, {\it i.e.},
\[
Y(1)_T = (Y(1)\backslash N(T), (\gamma\backslash\gamma_0)\cup\mu_T)\cong (M,\gamma),
\]
where $\mu_T$ is a meridian of $T$, and the assumption of Case 2.2 implies that $[T] = 0 \in H_1(Y\backslash B^3,\partial (Y\backslash B^3);\mathbb{Q})$. Then \cite[Proposition 1.4]{LY2020} concludes that
\[
\dim I^{\sharp}(Y) = \dim SHI(Y(1)) \leq \dim SHI(M,\ga)=2.
\]
To show (3), we need to unpack the proof of \cite[Proposition 1.4]{LY2020}, which is ultimately the proof of \cite[Proposition 3.14]{LY2020}. We view \cite[Proposition 1.4]{LY2020} as a special case of \cite[Proposition 3.14]{LY2020} when $T_0 = \emptyset$. \cite[Equation (3.2)]{LY2020} implies that there are sutures $\Gamma_{n-1}$ and $\Gamma_n$ such that we have an exact triangle
\begin{equation}\label{eq: exact triangle, 1}
	\xymatrix{
SHI(-M,-\Gamma_{n-1})\ar[rr]&&SHI(-M,-\Gamma_n)\ar[ld]\\
&SHI(-M,-\gamma)\ar[lu]&
}
\end{equation}
And \cite[Lemma 3.21]{LY2020} can be re-written (by replacing $n$ in the original equation by $n-1$) as
\begin{equation}\label{eq: exact triangle, 2}
	\xymatrix{
SHI(-M,-\Gamma_{n-1})\ar[rr]&&SHI(-M,-\Gamma_n)\ar[ld]\\
&I^{\sharp}(-Y)\cong SHI(-Y(1))\ar[lu]&
}
\end{equation}
Hence some basic linear algebra together with \eqref{eq: exact triangle, 1} and \eqref{eq: exact triangle, 2} implies that
\[
\dim I^{\sharp}(-Y)\equiv \dim SHI(-M,-\ga)~{\rm mod}~2
\]
As in \cite[Theorem 1.2]{li2018gluing}, we know $SHI(M,\gamma)$ and $SHI(-M,\gamma)$ are naturally dual to each other. Since $\partial M \cong T^2$, we know $\gamma$ and $-\gamma$ are isotopic, we conclude that
\[
\dim SHI(-M,-\ga) = \dim SHI (M,\gamma).
\]
Similar argument applies to $SHI(-Y(1))$ and we are done.
\epf

Recall $M = S^3\backslash N(J)$ is a knot complement and $Y$ is obtained from $M$ by filling along a component of $\gamma$, and hence $Y$ can be viewed as a Dehn surgery along $J$. The assumption of Case 2.2 implies that the surgery slope is non-zero. By passing to the mirror of $J$, which corresponds to reversing the orientation of $M$, we can assume that the surgery slope is positive. By Claim 4, we know that
$$\dim I^{\sharp}(Y)=2=|H_1(Y)|.$$
Note that by \cite[Theorem 1.15]{baldwin2019lspace}, the unknot and the right-handed trefoil are the only two knots on which the positive Dehn surgeries induce instanton L-spaces $Y$ with $|H_1(Y)|=2$. (According to the theorem, such a knot must be fibered and has genus at most $1$ and thus must be either the unknot, the trefoil, or the figure eight. Note the last knot is not strongly quasi-positive.) The case of unknot still reduces to Case 1. The case of the right-handed trefoil is a new one. By \cite[Theorem 1.1, Table 1]{baldwin2020concordance}, the surgery slope must be $2$. Hence $\ga$ consists of curves of slope $2$, which concludes the proof of the necessary condition.

Finally, the sufficient condition follows immediately from the first isomorphism in (\ref{eq: top grading and Seifert surface complement}) and the fact that gluing a product sutured manifold other than a $3$-ball to an arbitrary sutured manifold via identification of a suture does not change the sutured instanton Floer homology ({\it cf.} \cite[Proposition 6.7]{kronheimer2010knots}).
\epf

\brem\label{rem}
We have the following comments which strengthen the description of the guts in Theorem \ref{thm: rank-2 top grading}.
\begin{enumerate}
	\item We can compute the Euler characteristic in each of the three models. From \cite{kronheimer2010instanton}, for the first model, we have$$\chi(KHI(S^3,K,g(K)))=\chi(SHI(M,\ga))=0.$$As a result, we know that the symmetrized Alexander polynomial $\Delta_K(t)$ of $K$ has degree at most $g(K)-1$. On the other hand, if $(M,\ga)$ is one of the other two models, we can compute  as in \cite{LY2021} that
	$$\chi(SHI(M,\ga))=\pm2.$$
	As a result, we know that $\Delta_K(t)$ has degree $g(K)$ and the top non-zero coefficients are $\pm 2$.
	\item Let $S$ be a minimal genus Seifert surface of the knot $K\subset S^3$. Recall as in \eqref{eq: decompose along Seifert surface} and \eqref{eq: decompose along product annuli}, we have a decomposition
\[
S^3\backslash N(K)\stackrel{S}{\leadsto}(M^\p,\ga^\p)\stackrel{A}{\leadsto}(M,\ga)\sqcup (M_1,\ga_1)
\]
where $(M_1,\gamma_1)$ is a product sutured manifold and $(M,\gamma)$ is the guts. We write $(M_1,\gamma_1) = ([-1,1]\times F,\{0\}\times\partial F)$. The proof of \cite[Lemma 3.4]{baldwin2022nearly} implies that the Seifert surface complement $(S^3\backslash [-1,1]\times S,\{0\}\times\partial S)$ admits no product annuli whose boundary has a component that is parallel to the suture $\{0\}\times\partial S$ on $\partial (S^3\backslash [-1,1]\times S)$. As a result, we can further conclude that $\partial F$ must have one more component than $\ga$, and all but one components of $\partial F$ are glued to all of $\ga$. This actually rules out one model in the case $g(K)=1$ as in the following example.
\end{enumerate}
\erem

\begin{exmp}
	We keep the notations as in Remark \ref{rem}. When $g(K)=1$, we know that
	$$S\cong(R_{+}(\ga)\cup \{1\}\times F).$$
	Since in all three models we have $\chi(R_{+}(\ga))=0$, we know that
	$$\chi(F)=-1.$$
	From part (2) of the Remark \ref{rem}, we know that $\partial F$ has one more component than $\gamma$. Then $\chi(F)$ rules out the model in which $\ga$ has four components. As a result, we only have two models:
	\begin{itemize}
		\item $M$ is the complement of the unknot and $\ga$ consists of two curves of slope $2$.
		\item $M$ is the complement of the right handed trefoil and $\ga$ consists of two curves of slope $2$.
	\end{itemize} 
	Furthermore, in this case, the surface $F$ must be a pair of pants. Yet gluing such a thickened pair of pants to $(M,\ga)$ along two of the three boundary components is equivalent to gluing a product $1$-handle to $(M,\ga)$. Turning this around, we know that $(M,\ga)$ being one of the above two models is obtained from the complement of the Seifert surface by a disk decomposition. This coincides with the discussion in \cite[Section 1.2]{baldwin2022nearly} right above \cite[Theorem 5.1]{baldwin2022nearly}. Note that these two models do exist: for example, they give rise to the knot $5_2$ in Rolfsen's table and the 2-twisted Whitehead double of the right-handed trefoil with positive clasp.
\end{exmp}

%————End from here————
%\newpage
\bibliographystyle{alpha}
\bibliography{Seifert_accepted.bbl}

\end{document}

%% file: command.tex
\usepackage{amsmath}
\usepackage{amssymb}
\usepackage{extarrows}
\usepackage{mathabx}
\usepackage[all]{xy}%commutative graph
\usepackage{amsbsy}
\usepackage{graphicx}
\usepackage{appendix}
\usepackage{float}
\usepackage{subfigure}
\usepackage[numbers,sort&compress]{natbib}
\usepackage{graphicx}%pictures
\usepackage{amsthm}
\usepackage{geometry}
\usepackage{fancyhdr}%Heads
\usepackage{color} %colors
\usepackage{overpic}%pictures with text
\usepackage{mathabx}
\usepackage{tikz-cd}
\usepackage{enumerate}
\usepackage{mathrsfs}
\usepackage{colonequals}
\usepackage{microtype}
\usepackage{hyperref}
\usepackage{diagbox}

\newtheorem{thm}{Theorem}[section]
\newtheorem{cor}[thm]{Corollary}
\newtheorem{prop}[thm]{Proposition}
\newtheorem{lem}[thm]{Lemma}

\theoremstyle{definition}
\newtheorem{defn}[thm]{Definition}
\newtheorem{cons}[thm]{Construction}
\newtheorem{exmp}[thm]{Example}

\newtheorem*{fact}{Fact}

\newtheorem*{claim}{Claim}

\theoremstyle{remark}
\newtheorem{rem}[thm]{Remark}

\numberwithin{equation}{section}

\newcommand{\beq}{\begin{equation*}\begin{aligned}}
\newcommand{\eeq}{\end{aligned}\end{equation*}}
\newcommand{\bpf}{\begin{proof}}
\newcommand{\epf}{\end{proof}}
\newcommand{\bthm}{\begin{thm}}
\newcommand{\ethm}{\end{thm}}
\newcommand{\bprop}{\begin{prop}}
\newcommand{\eprop}{\end{prop}}
\newcommand{\bcor}{\begin{cor}}
\newcommand{\ecor}{\end{cor}}
\newcommand{\blem}{\begin{lem}}
\newcommand{\elem}{\end{lem}}
\newcommand{\bdefn}{\begin{defn}}
\newcommand{\edefn}{\end{defn}}
\newcommand{\bcons}{\begin{cons}}
\newcommand{\econs}{\end{cons}}
\newcommand{\bexmp}{\begin{exmp}}
\newcommand{\eexmp}{\end{exmp}}
\newcommand{\brem}{\begin{rem}}
\newcommand{\erem}{\end{rem}}
\newcommand{\bfa}{\begin{fact}}
\newcommand{\efa}{\end{fact}}
\newcommand{\benu}{\begin{enumerate}[(1)]}
\newcommand{\eenu}{\end{enumerate}}

\newcommand{\bdia}{\begin{displaymath}\xymatrix}
\newcommand{\edia}{\end{displaymath}}

\newcommand{\ga}{\gamma}

\newcommand{\p}{\prime}

\newcommand{\intg}{\mathbb{Z}}

\newcommand{\ra}{\rightarrow}

%% file: Seifert_accepted.bbl
\begin{thebibliography}{KM10b}

\bibitem[AZ22]{Agol2022guts}
Ian Agol and Yue Zhang.
\newblock {Guts in Sutured Decompositions and the Thurston Norm}.
\newblock {\em ArXiv:2203.12095, v1}, 2022.

\bibitem[BS21]{baldwin2020concordance}
John~A. Baldwin and Steven Sivek.
\newblock Framed instanton homology and concordance.
\newblock {\em J. Topol.}, 14(4):1113--1175, 2021.

\bibitem[BS22a]{baldwin2022nearly}
John~A. Baldwin and Steven Sivek.
\newblock Floer homology and non-fibred knot detection.
\newblock {\em ArXiv: 2208.03307, v1}, 2022.

\bibitem[BS22b]{baldwin2019lspace}
John~A. Baldwin and Steven Sivek.
\newblock Instanton and {L}-space surgeries.
\newblock {\em J. Eur. Math. Soc.}, 25(10):4033--4122, 2022.

\bibitem[Ghi08]{ghiggini2008knot}
Paolo Ghiggini.
\newblock Knot {F}loer homology detects genus-one fibred knots.
\newblock {\em Amer. J. Math.}, 130(5):1151--1169, 2008.

\bibitem[GL23]{li2019decomposition}
Sudipta Ghosh and Zhenkun Li.
\newblock Decomposing sutured monopole and instanton {F}loer homologies.
\newblock {\em Selecta Math. (N.S.)}, 29(3):Paper No. 40, 60, 2023.

\bibitem[Juh10]{juhasz2010polytope}
Andr\'{a}s Juh\'{a}sz.
\newblock The sutured {F}loer homology polytope.
\newblock {\em Geom. Topol.}, 14(3):1303--1354, 2010.

\bibitem[KM10a]{kronheimer2010instanton}
Peter~B. Kronheimer and Tomasz~S. Mrowka.
\newblock Instanton {F}loer homology and the {A}lexander polynomial.
\newblock {\em Algebr. Geom. Topol.}, 10(3):1715--1738, 2010.

\bibitem[KM10b]{kronheimer2010knots}
Peter~B. Kronheimer and Tomasz~S. Mrowka.
\newblock Knots, sutures, and excision.
\newblock {\em J. Differ. Geom.}, 84(2):301--364, 2010.

\bibitem[Li21a]{li2018gluing}
Zhenkun Li.
\newblock Gluing maps and cobordism maps in sutured monopole and instanton
  {F}loer theories.
\newblock {\em Algebr. Geom. Topol.}, 21(6):3019--3071, 2021.

\bibitem[Li21b]{li2019direct}
Zhenkun Li.
\newblock Knot homologies in monopole and instanton theories via sutures.
\newblock {\em J. Symplectic Geom.}, 19(6):1339--1420, 2021.

\bibitem[LY22]{LY2020}
Zhenkun Li and Fan Ye.
\newblock {Instanton Floer homology, sutures, and {H}eegaard diagrams}.
\newblock {\em J. Topol.}, 15(1):39--107, 2022.

\bibitem[LY23]{LY2021}
Zhenkun Li and Fan Ye.
\newblock Instanton {F}loer homology, sutures, and {E}uler characteristics.
\newblock {\em Quantum Topol.}, 14(2):201--284, 2023.

\bibitem[Ni07]{ni2007knot}
Yi~Ni.
\newblock Knot {F}loer homology detects fibred knots.
\newblock {\em Invent. Math.}, 170(3):577--608, 2007.

\bibitem[Sca15]{scaduto2015instanton}
Christopher Scaduto.
\newblock Instantons and odd {K}hovanov homology.
\newblock {\em J. Topol.}, 8(3):744--810, 2015.

\bibitem[XZ19]{xie2019tangle}
Yi~Xie and Boyu Zhang.
\newblock Instanton {F}loer homology for sutured manifolds with tangles.
\newblock {\em ArXiv:1907.00547, v2}, 2019.

\end{thebibliography}
